\documentclass[aip,cha,amsmath,amssymb,reprint,longbibliography,floatfix]{revtex4-1}
\usepackage{caption}
\usepackage{placeins}
\usepackage{graphicx}% Include figure files
\usepackage{dcolumn}% Align table columns on decimal point
\usepackage{bm}% bold math
\usepackage[a4paper,margin=2.5cm]{geometry}
\usepackage[utf8]{inputenc}
\usepackage[T1]{fontenc}
\usepackage{mathptmx}
\usepackage{etoolbox}
\usepackage{amsfonts}
\usepackage{amssymb}
\usepackage{graphicx,color}
\usepackage{amsmath}
\usepackage{amsthm}
\usepackage{subcaption}
\graphicspath{{./Figs/}}
\DeclareGraphicsExtensions{.pdf,.jpg,.png,.eps,.gif}
\def\ds{\displaystyle}
\makeatletter
\def\@email#1#2{%
 \endgroup
 \patchcmd{\titleblock@produce}
  {\frontmatter@RRAPformat}
  {\frontmatter@RRAPformat{\produce@RRAP{*#1\href{mailto:#2}{#2}}}\frontmatter@RRAPformat}
  {}{}
}%
\newtheorem{theorem}{Theorem}[section]

\newtheorem{proposition}[theorem]{Proposition}

\theoremstyle{definition}

\makeatother
\begin{document}

\preprint{AIP/123-QED}

\title[Legendre Functions and the Non-Integrability of a Hamiltonian System]{\bf  Legendre Functions and the Non-Integrability of a Hamiltonian System}
% Force line breaks with \\
\author{D. Neykova}
 \email{ dneykova@uni-sofia.bg}
 \thanks{ORCID: https://orcid.org/0009-0008-7892-4031}
\author{G. Georgiev}%
 \email[Author to whom correspondence should be addressed:]{ggeorgiev3@fmi.uni-sofia.bg}
 \thanks{ORCID: https://orcid.org/0000-0002-8918-1420}
\affiliation{ 
 Faculty of Mathematics and Informatics \\
 Sofia University ``St. Kliment Ohridski'', \\
  5 James Bourchier Blvd., 1164 Sofia, Bulgaria
}
%

%\date{\today}% It is always \today, today,
             %  but any date may be explicitly specified

\begin{abstract}
We study the meromorphic integrability of a two-dimensional Hamiltonian system with a homogeneous polynomial potential of degree six. Our approach is based on the Ziglin–Morales–Ruiz–Ramis–Simó theory, which relates integrability of Hamiltonian systems to the differential Galois group of the variational equations along a non-equilibrium particular solution. For Hamiltonian systems with homogeneous potentials, the first variational equations are typically reduced to hypergeometric equations; however, the standard approach requires the computation of Darboux points, which may be difficult in practice. In the present work this step is avoided by reducing the variational equation to an associated Legendre equation. We then use known criteria for the differential Galois group of the associated Legendre equation to obtain non-integrability conditions for a Hamiltonian system with a homogeneous polynomial potential of degree six. Special attention is paid to the parameter regimes not fully resolved by the classical criteria. In the solvable cases of the first variational equation, second variational equations are used to identify conditions under which a non-zero logarithmic term appears in the local solutions. This shows that the identity component of the corresponding differential Galois group is not Abelian, and hence the Hamiltonian system admits no additional meromorphic first integral.
\end{abstract}

\maketitle

\begin{quotation}
Hamiltonian systems with homogeneous polynomial potentials of degree six arise in the modeling of a variety of complex nonlinear phenomena, including molecular vibrations, nonlinear electrical circuits, and approximate cosmological dynamics. The study of chaotic behavior in such systems constitutes a natural extension of the classical analysis of lower-degree polynomial potentials. While quadratic potentials correspond to integrable systems, fourth-degree potentials exhibit both integrable and non-integrable dynamics, thereby marking the onset of complex chaotic behavior. In contrast, sixth-degree polynomial potentials are characterized predominantly by non-integrability and the widespread emergence of chaos, making them an important subject in the qualitative theory of dynamical systems and nonlinear mechanics.
\end{quotation}

\section{Introduction}
We consider Hamiltonian systems of type
\begin{equation}
\label{HomPot}
 H=\sum_i\frac{1}{2}p_i^2+V(q),
\end{equation}
where $q=(q_1,\dots , q_n)$ and $p=(p_1,\dots , p_n)$ are canonical coordinates in  a symplectic linear space $\mathbb{C}^{2n}$ and $V(q)$  is a homogeneous polynomial function. If $F=F(q,\,p)$ and $G=G(q,\,p)$ are two functions, then their Poisson bracket is defined by $\{F,\,G\}=\sum_i(\frac{\partial F}{\partial q_i}\frac{\partial G}{\partial p_i}-\frac{\partial F}{\partial p_i}\frac{\partial G}{\partial q_i})$. The functions  $F$  and $G$ are in involution if $\{F,\,G\}=0$ is satisfied. Non-constant function $F$
 is the first integral for the Hamiltonian system with Hamiltonian $H$ if $\{F,\,H\}=0$.  Since the Poisson bracket is anti-symmetric, it is clear that $H$ is always the first integral. 
A Hamiltonian system is called integrable if a set of $n$ independent first integrals in involution exists for it (see Appendix A for further details). 

The study of Hamiltonian systems with a homogeneous polynomial potential is an important problem in modern mechanics and physics. In our opinion, the reason for this research interest is the quite important physical problems can be considered as Hamiltonian systems with a homogeneous potentials. The integrability of Hamiltonian systems of type (\ref{HomPot}) has been studied for a long time, but the results of this work are limited to a rather small class of solved problems. The approaches to the study of this problem are: direct method \cite {Yarmo1}, \cite {Yarmo2} and Painleve test \cite{GRAMMATICOS1}, \cite{GRAMMATICOS2}, \cite{GRAMMATICOS3}, \cite{GRAMMATICOS4}, \cite{GRAMMATICOS5} and \cite{GRAMMATICOS6}. The most remarkable and perhaps most general result was obtained in \cite{Yoshida0}. This paper investigates the necessary conditions for integrability in terms of analytic first integrals.
In 1982 in \cite{Ziglin1} and 1983 in \cite{Ziglin2}, Ziglin reduced the study of the integrability of Hamiltonian systems to the investigation of the monodromy group of the linear equations obtained by the variation of the system in a neighborhood of its non-equilibrium solution. The Ziglin's theory gives the opportunity for the study complex Hamiltonian systems for integrability. Using this theory, Yoshida \cite{Yoshida1} was able to obtain a criterion for the integrability of two-dimensional systems with a homogeneous polynomial potential (the necessary conditions  in terms of holomorphic first integrals). Later in \cite{Yoshida2} Yoshida generalized his result to $n$ degrees of freedom.
In the late twentieth century, Ziglin's theory was linked to the differential Galois theory in the works of Churchill, Baider, Rod, Singer  in \cite{Church1}, \cite{Church2}, Morales-Ruiz, Ramis, and Simo in \cite{MR1}, \cite{MR2}, \cite{MRS1}. 
The idea here is to move from studying the discrete monodromy group to the continuous algebraic  differential Galois group for linear variational equations. The closure of the monodromy group in the Zariski topology gives us the differential Galois  group, as noted in the Picard--Vessiot theory (see Appendix A.)

The use of differential Galois theory to determine the integrability of Hamiltonian systems with homogeneous potentials was studied in\cite{MR1}. In this publication, necessary conditions for integrability in terms of meromorphic first integrals were formulated. In the papers \cite{MPr1} and \cite{MPr2}, the cases of homogeneous polynomial potentials of degrees 3 and 4 are investigated. Homogeneous potentials with a negative degree of homogeneity are considered in \cite{Llibre} and \cite{Valls}. 

The approach implemented to study the non-integrability of $V_6$ in this paper differs from \cite{MR1}, \cite{MPr1} and \cite{MPr2}. Here, the search for Darboux points is skipped. With a change of variables $VE_1$, which is a Fuchsian differential equation with 5 regular singularities, is reduced to a Legendre equation with 3 singularities.

We examine systems with Hamiltonian
$H=\sum_i\frac{1}{2}p_i^2+V(q)$, with potential $V(q)$ shown in the form $V(q)=\sum_l V_l(q)$, here $V_l (q)$ are homogeneous polynomials of degree $l\in \mathbb{N}$, $l\ge 2$.

Every homogeneous polynomial potential $V(q)$  is represented as a sum $V(q)=V_{min}(q)+\dots +V_{max}(q)$ of homogeneous components of increasing degrees. We  call the potential $V(q)$ integrable if the Hamiltonian system $H=\sum_i\frac{1}{2}p_i^2+V(q)$ is integrable. If a homogeneous potential $V(q)=V_{min}(q)+\dots +V_{max}(q)$ is integrable then $V_{min}(q)$ and $V_{max}(q)$ are also integrable. This remarkable property can be seen in \cite{Yarmo2}.

In this paper, we investigate the meromorphic non-integrability of a sixth-degree homogeneous polynomial potential, aiming to extend the previously known results.
We assume that $\mathbb{C}^4$ is a symplectic linear space with canonical variables $q=(q_1,q_2)$ and $p=(p_1,p_2)$. We  study the Hamiltonian systems defined by the Hamiltonian $H=\frac{1}{2}\sum_1^2p_i^2+V_l(q)$, where $V_l(q)$ is a homogeneous function of degree $l$. The equations of motion of the above system are
\begin{eqnarray*}
\dot{q}_i & = & p_i,\\
\dot{p}_i &=& -\frac{\partial V(q_1,q_2)}{\partial q_i} ,\, i=1,\, 2.           
\end{eqnarray*}
We say that a Hamiltonian system with 2 degrees of freedom is integrable in the sense of Liouville if there exists a second (other than $H$) first integral  $F$ for which the gradients  $dH$ and $dF$  are linearly independent (see \cite{Llibre} for details).
We study two-dimensional Hamiltonian systems for the absence of an additional meromorphic first integral of motion.
Let us consider two dimensional model with sixth-order homogeneous potential
\begin{eqnarray}
\label{1.1}
H & = &\frac{1}{2}(p_r^2+p_z^2) \nonumber\\
& + & r^6+Ar^4z^2+Dr^3z^3+Br^2z^4+Erz^5+Cz^6,
\end{eqnarray}
 where $A$, $B$, $C$, $D$ and $E$ are appropriate
real constants and $(r,z,p_r,p_z)$ are in $\mathbb{C}^4$, for existing  additional  meromorphic integral.

Let us justify the absence of the coefficient in front of $r^5z$: If we assume that this coefficient in $V_6$ is $\alpha $, i.e. we have
$$V_6=r^6+\alpha r^5z+Ar^4z^2+Dr^3z^3+Br^2z^4+Erz^5+Cz^6.$$
Now let us change the variables $r\rightarrow \frac{6r-\alpha z}{36}$ and $z\rightarrow\frac{z}{6}$ (a linear change,  we can also define a canonical one of the Hamiltonian system). This is enough to change the coefficient $\alpha$  to 0 (see \cite{Yarmo1}). 

The Hamiltonian equations are:
\begin{align}
\label{1.2}
\dot r  &=  p_r,\nonumber\\
\dot p_r & =  -(6r^5+4Ar^3z^2+3Dr^2z^3+2Brz^4+Ez^5),\nonumber  \\
\dot z  &=  p_z,\nonumber\\
\dot p_z  & =  -(2Ar^4z+3Dr^3z^2+4Br^2z^3+5Erz^4+6Cz^5),
\end{align}
 (here as usual $\dot{}=\frac{d}{dt}$).

It is convenient to use the following notation:
$$2\mathbb{Z}_{\ge 0}:=\{0,\, 2, \, 4, \, 6 \dots \},$$
$$2\mathbb{Z}_{\ge0}+1:=\{1,\, 3, \, 5, \, 7 \dots \},$$
$$3\mathbb{Z}_{\ge0}:=\{0,\, 3, \, 6, \, 9 \dots \}.$$
The main result in this work is:
\begin{theorem}
\label{Th_Main}
Let $\tau=A$, or $\tau=\frac{B}{C}$, ($C\ne 0$) (here "or" is exclusive). Then, for $\tau\ne 3$ the system (\ref{1.2}) does not have an additional meromorphic first integral, if at least one of the following conditions holds:

1. $\sqrt{\tau+1}\notin\mathbb{Q}$;

 2. $\sqrt{\tau+1}\in\mathbb{Q}\setminus\{ 3k+1, -3k+2, 6k+3\}$,  for $k\in\mathbb{Z}$;
 
 3. $\sqrt{\tau+1}\in\{3k+1, -3k+2, 6k+3\}$, for $k\in\mathbb{Z}$;

3.1. if  $\sqrt{\tau+1}=3k+1$,  the conditions are:  $k\in 2\mathbb{Z}_{\ge 0}$, and $D\ne 0$;

3.2. if  $\sqrt{\tau+1}=-3k+2$, then $k\in 2\mathbb{Z}_{\ge0}+1$, and $D\ne 0$;

3.3. if $\sqrt{\tau+1}=6k+3$, then $k\in 3\mathbb{Z}_{\ge0}$.

\end{theorem}

This study may be viewed as a complement to the classical results presented in \cite{Yoshida1} and \cite{Yoshida2}. Although some questions remain open and new ones arise, as discussed in the concluding section, the present paper introduces a novel approach for investigating parameter regimes that are not fully resolved by the existing theory. In particular, for a broad range of parameter values, we analyze cases in which the variational equations are solvable.

The paper is organized as follows. In Section II, we investigate the conditions for the absence of Liouvillian solutions of the Legendre equation. In Section III, we derive a suitable particular solution of the system and determine the corresponding variational equations $VE_1$ in a neighborhood of this solution. We then relate these equations to the Legendre equation. In Section IV, we analyze the cases that are not covered by the results of Section III and prove the main theorem of the paper. This is followed by a discussion of the results obtained and illustrative examples for $\tau = 3$. Finally, Appendix A provides the necessary background on Differential Galois Theory.

\section{The associated Legendre equation}

The object of study of differential Galois theory (Picard--Vessiot theory) is the integrability of linear differential equations and systems. A linear differential equation is said to be integrable if its general solution can be expressed in terms of algebraic functions, exponentials of quadratures, and quadratures over a suitable field of constants. The integrability properties of such equations are encoded in the identity component of their differential Galois group: the equation is integrable if and only if the identity component of its Galois group is solvable.

Let us now consider an integrable Hamiltonian system $X_H$ and the variational equations obtained by linearizing the system in a neighborhood of a non-equilibrium particular solution. The resulting linear system, denoted by $VE$, must also be integrable in the sense of differential Galois theory. Equivalently, the corresponding Picard--Vessiot extension must be Liouvillian, or, in algebraic terms, the identity component of the associated differential Galois group must be a solvable algebraic group.

A stronger statement is provided by the Ziglin--Morales--Ramis theory: in the case of complete integrability, the identity component of the differential Galois group of $VE$ is necessarily Abelian. In other words, the Abelian structure of the Poisson algebra of first integrals of $X_H$ is reflected at the level of the variational equations.
The connection between these two different concepts of integrability is the theory of Ziglin-Morales-Ruiz-Ramis-Simo.  (See \cite{MR1} and Appendix A for details.)

In this text, a solvable group will be understood as an algebraic group represented by upper triangular matrices.

In this section, we study the absence of Liouvillian solutions (solutions that are expressed through elementary functions and their integrals) to the associated Legendre equation
\begin{eqnarray}
\label{Legendre}
(1-z^2)\frac{d^2 w}{dz^2}-2z\frac{d w}{dz}+\left(p(p+1)-\frac{q^2}{1-z^2}\right)w=0,\\ \nonumber z\in\mathbb{C},\, p, \,q\in \mathbb{R},\, p\pm q\notin \mathbb{Z}_{\leq 0} ,
\end{eqnarray}
and we apply the result to study the potential
 $$V_6(r,z)=r^6+Ar^4z^2+Dr^3z^3+Br^2z^4+Erz^5+Cz^6,$$
  ($A,\, B,\, C,\, D ,\, E \in \mathbb{R}$) for meromorphic non-integrability (within the scope of Liouville).

The restriction $p\pm q\notin \mathbb{Z}_{\leq 0}$ is necessary because, for these parameter values, the solutions of (\ref{Legendre}) become linearly dependent over $\mathbb{C}$, as shown below.

  The study plan is as follows: 
  
First, we determine the exponents at the singularities $\pm 1$ and $\infty$, and use them to derive the generators of the monodromy group. This enables us to identify the cases in which the monodromy group is non-commutative.

Second, by expressing the solutions in terms of the hypergeometric function, we apply the solvability criteria of \cite{Kimura} and \cite{Primitivo}.

In the subsequent sections, we investigate the solvable but non-commutative cases.
 
Let us write down the known facts about the solutions of (\ref{Legendre}) and outline some properties required for our purposes. (We follow \cite{DLMF}.)
With
\begin{eqnarray}
\label{solsLeg}
P_p^q(z) & = & \left(\frac{1+z}{1-z}\right)^{q/2} \,  {_2}F_1\, (p+1,-p;1-q,\frac{1}{2}-\frac{z}{2}), \nonumber\\
Q_p^q(z) & = & \frac{\pi}{2\sin(q\pi)}(\cos(q\pi)\left(\frac{1+z}{1-z}\right)^{q/2} \nonumber\\
& \times & {_2}F_1\, (p+1,-p;1-q,\frac{1}{2}-\frac{z}{2}) \nonumber \\
&- &\left(\frac{1-z}{1+z}\right)^{q/2}\frac{\Gamma(p+q+1)}{\Gamma(p-q+1)} \nonumber\\
& \times & {_2}F_1\, (p+1,-p;1+q,\frac{1}{2}-\frac{z}{2}))    
\end{eqnarray}
we can show the solutions of (\ref{Legendre}) ($P$ and $Q$ Legendre functions), shown by using the hypergeometric function $_2F_1\, (a,b;c,z)$.  
The singularities of equation (\ref{Legendre}) are the points $z=-1$, $z=1$ and $z=\infty$, which are regular. The following equalities are valid

\begin{eqnarray}
\label{MonodEq}
P_p^{-q}(ze^{si\pi }) & = & e^{s p i\pi }P_p^{-q}(z)+\frac{2i\sin{(p+1/2)s\pi}}{\cos{p\pi}\Gamma(q-p)}Q_p^q(z)\nonumber\\
Q_p^q(ze^{si\pi }) & = &(-1)^s e^{-s p i\pi }Q_p^q(z),\,s\in \mathbb{Z}.
\end{eqnarray}
Let us focus on the generators of the monodromy group of equation (\ref{Legendre}). 
For this purpose, we are looking for the indicative equations for the singular points $\pm 1$ and $\infty$. 

We have  $\rho^2-\rho +\ds{\frac{1-q^2}{4}}=0$ (see \cite{Poole}) with a roots (exponents) $\rho_{1,2}=\ds{\frac{1\pm q}{2}}$,
 and for $\infty$ it is $\lambda^2+\lambda-p^2-p=0$ with an exponents $\lambda_1=-p-1$, $\lambda_2=p$. 
 
Therefore, we can specify the generators of the monodromy group. For the points $\pm 1$ and $\infty$ these are the matrices 
 $$M_{-1}=\begin{pmatrix}
1 & 0  \\
2(\cos2\pi p-\cos2\pi q) & {e^{-2i\pi q}}
\end{pmatrix},$$  
$$M_{1}=\begin{pmatrix}
1 & {e^{-2i\pi q}}-1  \\
0 & {e^{-2i\pi q}}
\end{pmatrix},$$ 
and 
$$M_{\infty}=\begin{pmatrix}
{e^{-2i\pi q}}(2\cos 2\pi p -1) & 0  \\
{e^{-2i\pi q}}(2\cos 2\pi p -1)-1& e^{-2\pi i q }
\end{pmatrix}.$$

It is known from the theory (see Appendix A) that a necessary condition for the existence of Liouvillian solutions of a linear differential equation is that its monodromy group be commutative. Using the generators obtained above, we conclude that equation (\ref{Legendre}) does not admit Liouvillian solutions whenever at least one of the parameters $p$ or $q$ is irrational.
\begin{proposition}
\label{NonRational}
The  equation (\ref{Legendre}) does not have Liouvillian solutions if
 at least one of $p$,  $q\notin\mathbb{Q}$.
\end{proposition}

In the remainder of this paper, we assume that $p,q\in\mathbb{Q}$. 
We now have explicit expressions for the solutions of (\ref{solsLeg}) in terms of the Legendre functions $P_p^q(z)$ and $Q_p^q(z)$, which can also be represented by the hypergeometric function
\[
{}_2F_1\left(p+1,-p;1-q,\frac{1}{2}-\frac{z}{2}\right).
\]
This allows us to apply the solvability criteria of Kimura \cite{Kimura} and \cite{Primitivo}. Thus we get the following
\begin{theorem}
\label{Th_Leg}
({\bf \cite{Primitivo} Acosta-Humanez et al.)} 
Let $p,q\in\mathbb{Q}$. Then equation (\ref{Legendre}) admits Liouvillian solutions if at least one of the following conditions holds:

\begin{enumerate}
\item At least one of the numbers $q+p$, $q-p$, or $p$ is an integer;

\item The quantities $\pm q$, $\pm p$, and $\pm (p+q)$ belong to one of the following families:

\begin{enumerate}
\item[(a)] $\mathbb{Z}+\frac{1}{2}$, $\mathbb{C}$, $\mathbb{C}$;

\item[(b)] $\mathbb{Z}\pm\frac{1}{3}$, $\frac{1}{2}\mathbb{Z}\pm\frac{1}{3}$, $\mathbb{Z}+\frac{1}{6}$;

\item[(c)] $\mathbb{Z}\pm\frac{2}{5}$, $\frac{1}{2}\mathbb{Z}\pm\frac{1}{5}$, $\mathbb{Z}+\frac{1}{10}$;

\item[(d)] $\mathbb{Z}\pm\frac{1}{3}$, $\frac{1}{2}\mathbb{Z}+\frac{2}{5}$, $\mathbb{Z}+\frac{1}{10}$;

\item[(e)] $\mathbb{Z}\pm\frac{1}{5}$, $\frac{1}{2}\mathbb{Z}\pm\frac{2}{5}$, $\mathbb{Z}+\frac{1}{10}$;

\item[(f)] $\mathbb{Z}\pm\frac{2}{5}$, $\frac{1}{2}\mathbb{Z}\pm\frac{1}{3}$, $\mathbb{Z}+\frac{1}{6}$.
\end{enumerate}
\end{enumerate}
\end{theorem}

It is not difficult to conclude that, if $q=0$, then the Legendre equation (\ref{Legendre}) does not admit Liouvillian solutions for $p\notin\mathbb{Z}$.

\section{Sixth-order homogeneous potential}

To investigate the meromorphic non-integrability of the Hamiltonian system (\ref{1.2}), we apply the Ziglin--Morales--Ruiz--Ramis--Sim\'o theory.

First, we find a partial solution of (\ref{1.2}). Setting
$z=p_z=0$ in (\ref{1.2}), we obtain
\begin{equation}
\label{1.33}
\ddot r=-6r^5,
\end{equation}
multiplied by $\dot r$ and integrated by $t$ we get
\begin{equation}
\label{1.3}
{\dot r}^2=-2(r^6+h^6),
\end{equation}
where $h$ is a constant. We should note that the functions $r(t)$ are finitely branched.

Next, we follow the Ziglin--Morales--Ramis approach and consider the invariant manifold
\[
(r,p_r,z,p_z)=(r(t),\dot r(t),0,0),
\]
where $r(t)$ satisfies (\ref{1.3}). It is convenient to choose the field of constants to be
$K = \mathbb{C}(r)$ - rational functions over a complex variable $r$. After deriving the variational equations $(VE)$, we obtain $\xi_{11}=dr$, $\eta_{11}=dp_r$,
$\xi_{12}=dz$,  $\eta_{12}=dp_z$,  and we achieve:
\begin{align}
\label{VE1.5}
 \ddot \xi_{11 }&=  -2Ar^4\xi_{11},\nonumber\\
\ddot \xi_{12}  &=  -30r^4\xi_{12}.
\end{align}

Now we perform the change of variable
\[
t\mapsto r(t),
\]
for which
\[
\ddot{\xi}=\xi''\dot r^{\,2}+\xi'\ddot r.
\]
From (\ref{1.33}) and (\ref{1.3}), we derive the equations $VE_1$.. Here and in what follows, the notation
\[
'=\frac{d}{dr}
\]
is used. Thus, the equations $VE_1$ are reduced to the following two Fuchsian linear differential equations:
\begin{eqnarray}
\label{VE1.6}
 \xi_{11}''  +\frac{3r^5}{(r^6+h^6)}\xi_{11}'  -\frac{Ar^4}{(r^6+h^6)}\xi_{11}=0, \nonumber \\
 \xi_{12}''  +\frac{3r^5}{(r^6+h^6)}\xi_{12}'
-\frac{15r^4}{(r^6+h^6)}\xi_{12}=0.
\end{eqnarray}
Equations ( \ref{VE1.6}) are a Fuchsian and have seven regular singularities - the roots of $r^6+h^6=0$ and $\infty$. 
There are three possible ways to investigate (\ref{VE1.6}) the solvability:
There are three possible approaches for investigating the solvability of (\ref{VE1.5}). The first is to apply the Kovacic algorithm \cite{BaBi} to the first equation. Although this approach is technically involved, it remains applicable. The second approach is to study the existence of hyperexponential solutions \cite{Bouchier1}, \cite{Bouchier2}. The third approach is to perform the change of variables
\begin{equation}
\label{ChangeVarVE}
y^2:=(1+\frac{r^6}{h^6}),\, {\Xi}_{11}:=\frac{\xi_{11}}{{(h^6(y^2-1))^{1/12}}},
\end{equation}
we get an easy to investigate equation. We obtain
\begin{equation}
\label{NVE1.6_New}
\frac{d^2\Xi_{11}}{dy^2}  -\frac{2z}{1-y^2}\frac{d\Xi_{11}}{dy} +\left( \frac{4A-5}{36}-\frac{1/36}{1-y^2}\right)\Xi_{11}=0.
\end{equation}
This change of variables defines a finitely branched covering
\[
\mathbb{CP}^1 \to \mathbb{CP}^1.
\]
In general, the differential Galois group changes under such a transformation, whereas its identity component remains unchanged (see \cite{MR1}, p.~28). In this way, the study of the original problem is reduced to a significantly simpler one. We first eliminate the cases for which (\ref{NVE1.6_New}) does not admit Liouvillian solutions and then focus on the remaining solvable cases.
The equation (\ref{NVE1.6_New}) is associated Legendre equation with $p=-\ds{\frac{1}{2}}+\ds{\frac{\sqrt{1+A}}{3}}$ and $q=\ds{\frac{1}{6}}$. We apply the results of Proposition \ref{NonRational} and Theorem \ref{Th_Leg}:
\begin{proposition}
\label{NonRationalNEW}
The  system  (\ref{1.2}) is non-integrable for
$\sqrt{1+A}\notin\mathbb{Q}$.

\end{proposition}

\begin{proposition}
\label{NonZNEW}
For $\sqrt{1+A}\in\mathbb{Q}$, the  system  (\ref{1.2}) is non-integrable for
$\sqrt{1+A}\notin\{3k+1,\,-3k+2,\, 6k+3\}$, $ k\in\mathbb{Z}$.
\end{proposition}

We emphasize that the solvability of the differential Galois group is not sufficient for the existence of an additional meromorphic first integral. In order to prove non-integrability, it is necessary to determine whether the identity component of the differential Galois group is non-Abelian.

The non-solvable cases are now completely characterized.
We next investigate the solvable cases
\begin{align*}
\sqrt{1+A} &= 3k+1,\\
\sqrt{1+A} &= -3k+2,\\
\sqrt{1+A} &= 6k+3,
\end{align*}
where $k\in\mathbb{Z}$,
in order to derive additional non-integrability conditions.

%%%%%%%%%%
\section{The solvable cases}

In this section we are considering cases $\sqrt{1+A}=3k+1$, $\sqrt{1+A}=-3k+2$, and $\sqrt{1+A}=6k+3$ for $ k\in\mathbb{Z}$ that complete the research from the previous section.

Let us find the second variations of the Hamiltonian system with Hamiltonian (\ref{1.1}). We denote it with
\begin{eqnarray*}
r &= & r(t)+\varepsilon\xi_{11}+\varepsilon^2 \xi_{21}+\dots , \\
z & = & \varepsilon\xi_{12}+\varepsilon^2 \xi_{22}+\dots , \\
 p_r & = & \dot{r}(t)+\varepsilon\eta_{11}+\varepsilon^2 \eta_{21}+\dots\\
p_z & = & \varepsilon\eta_{12}+\varepsilon^2 \eta_{22}+\dots, 
 \end{eqnarray*}
here $(r,\, p_r,\, z, \, p_z,)=(r(t),\, \dot{r}(t),\, 0, \,0)$ is an invariant manifold of the system (\ref{1.2}).

We do the replacement in the system (\ref{1.2}) and  we compare the coefficients at $\varepsilon^2$. After changing the variables $t\rightarrow r(t) $ we  obtain

\begin{eqnarray}
\label{VE2}
 \xi_{21}''  & + & \frac{3r^5}{(r^6+h^6)}\xi_{21}'  -\frac{Ar^4}{(r^6+h^6)}\xi_{21} =K_2^{(1)}\nonumber\\
 \xi_{22}''  & + & \frac{3r^5}{(r^6+h^6)}\xi_{22}'
-\frac{15r^4}{(r^6+h^6)}\xi_{22} = K_2^{(2)},
\end{eqnarray}
when
\begin{eqnarray*}
K_2^{(1)} & = & 4A\frac{r^3\xi_{11}\xi_{12}}{(r^6+h^6)}+\frac{3D}{2}\frac{r^3(\xi_{11})^2}{(r^6+h^6)},
\end{eqnarray*}

\begin{eqnarray*}
K_2^{(2)} & = & 2A\frac{r^3(\xi_{11})^2}{(r^6+h^6)}+30\frac{r^3(\xi_{12})^2}{(r^6+h^6)}.
\end{eqnarray*}

Since the exponents of the finite singular points of VE1 (the roots of $r^6+h^6=0$ ) are 0 and 1, we consider VE1 and VE2 at the  singular point $\infty$ of the Riemann sphere. To this end, we change variables $r=\frac{1}{x}$ into (\ref{VE1.6}) and (\ref{VE2}) (here we assume $'=\frac{d}{dx}$). It is convenient to fix the field of constants $\mathbb{C}(x)$ - the rational functions of $x$. Here $x$ is a complex variable.
We have
\begin{eqnarray}
\label{VE1_infty}
 \xi_{11}''  &+& \frac{2h^6x^6-1}{x(h^6x^6+1)}\xi_{11}'  -\frac{A\xi_{11}}{x^2(h^6x^6+1)}=0, \nonumber \\
 \xi_{12}''  &+&  \frac{2h^6x^6-1}{x(h^6x^6+1)}\xi_{12}'
-\frac{15\xi_{12}}{x^2(h^6x^6+1)}=0,
\end{eqnarray}
\begin{eqnarray}
\label{VE2_infty}
 \xi_{21}''  & + & \frac{2h^6x^6-1}{x(h^6x^6+1)}\xi_{21}'   -\frac{A\xi_{21}}{x^2(h^6x^6+1)} \nonumber \\
                & = & 4A\frac{\xi_{11}\xi_{12}}{x(h^6x^6+1)}+\frac{3D}{2}\frac{(\xi_{11})^2}{x(h^6x^6+1)},\nonumber \\
 \xi_{22}''  &+&  \frac{2h^6x^6-1}{x(h^6x^6+1)}\xi_{22}'
-15\frac{\xi_{22}}{x^2(h^3x^3+1)}\nonumber\\
& = & 2A\frac{(\xi_{11})^2}{x(h^6x^6+1)}+30\frac{(\xi_{12})^2}{x(h^6x^6+1)}.
\end{eqnarray}
For further purposes, it is convenient to use the notation:
\begin{eqnarray}
\label{K_xi}
{K_2}^{(1)}(\xi_{11},\xi_{12}) := 4A\frac{\xi_{11}\xi_{12}}{x(h^6x^6+1)}\nonumber\\
+\frac{3D}{2}\frac{(\xi_{11})^2}{x(h^6x^6+1)}\\
{K_2}^{(2)}(\xi_{11},\xi_{12}):= 2A\frac{(\xi_{11})^2}{x(h^6x^6+1)}\nonumber\\
+30\frac{(\xi_{12})^2}{x(h^6x^6+1)}.
\end{eqnarray}

We also need to give the solutions of (\ref{VE1_infty}) in series near $0$.
\begin{align}
\label{xi_Sol_0}
\xi_{11}^{(1)}(x)
&= x^{1-\sqrt{1+A}}
\Bigg(
1+
\frac{A+3-3\sqrt{A+1}}
     {12(\sqrt{A+1}-3)}
\, h^6 x^6
+\dots
\Bigg),
\nonumber\\[4pt]
\xi_{11}^{(2)}(x)
&= x^{1+\sqrt{1+A}}
\Bigg(
1-
\frac{A+3+3\sqrt{A+1}}
     {12(\sqrt{A+1}+3)}
\, h^6 x^6
+\dots
\Bigg).
\end{align}
\begin{eqnarray}
\xi_{12}^{(1)}(x)=x^5\left( 1+\dots\right)\nonumber,\\
\xi_{12}^{(2)}(x)=x^{-3}\left( -20331280+\dots\right).
\end{eqnarray}

We would like to underline some known facts we will need for our further research. Solutions should be provided for  the equations (\ref{VE1_infty}) or (\ref{VE2_infty})  where Wronsky determinant is a constant (we need $const=1$). This happens exactly when the coefficient in front of the first derivative  in the differential equation disappears. 
We are transforming the linear differential equation 
\begin{equation}
\label{homNonNorm}
\xi (x)''+a(x)\xi (x)'+b(x)\xi (x)=0
\end{equation}
into normal form (with change $\xi (x)=\zeta (x)e^{-\frac{1}{2}\int{a(x)dx}}$).
we obtain the linear differential equation 
\begin{equation}
\label{homNorm}
\zeta (x)''-r(x)\zeta=0,
\end{equation}
where $r(x)=\frac{1}{2}a(x)'+\frac{1}{4}(a(x))^2-b(x)$.
We do the same transformation for a non-homogeneous linear differential equation (the same equation, but with added right-hand side $K(x)$), and we have 
$$ \xi (x)''+a(x)\xi (x)'+b(x)\xi (x)=K_2(x),$$
and the corresponding transformed equation
$$\zeta (x)''-r(x)\zeta=K_2(x)e^{\frac{1}{2}\int{a(x)dx}}.$$
We could, without restriction, assume that the equation (\ref{homNonNorm}) has a regular singularity at the point $0$. Let us denote by $\zeta_i(x)$, $i=1,\,2, \,3$ three randomly chosen solutions of (\ref{homNorm}). 
It is clear that $\xi_i (x)=\zeta_i (x)e^{-\frac{1}{2}\int{a(x)dx}}$, ($i=1,\,2, \,3$) are a solutions of ({\ref{homNonNorm}).
Now let us express the product $\zeta_1 (x) \zeta_2(x)\zeta_3(x) e^{\frac{1}{2}\int{a(x)dx}}$, in terms of $\xi_1(x)\xi_2(x)\xi_3(x)$.
We have
\begin{eqnarray*} 
\zeta_1 \zeta_2\zeta_3 e^{\frac{1}{2}\int{a(x)dx}}& = & \zeta_1 e^{-\frac{1}{2}\int{a(x)dx}}\\
&\times &\zeta_2 e^{-\frac{1}{2}\int{a(x)dx}}\\
&\times &\zeta_3 e^{-\frac{1}{2}\int{a(x)dx}}e^{\frac{1}{2}\int{a(x)dx}}(e^{\frac{3}{2}\int{a(x)dx}})\\
& = & \xi_1(x)\xi_2(x)\xi_3(x)e^{2\int{a(x)dx}}.
\end{eqnarray*} 

Now we change the form of (\ref{VE1.6}) and (\ref{VE2}). Let us replace $\xi $ with $\zeta $ by
\begin{equation*}
 \xi =\zeta e^{-\frac{1}{2}\int{ \frac{2h^6x^6-1}{x(h^6x^6+1)}}dx}=\zeta .\frac{x^{1/2}}{(h^6x^6+1)^{\frac{1}{4}}},
 \end{equation*}
  and we obtain
 \begin{eqnarray}
\label{VE2_zeta}
\zeta_{11}''   -   r_1(x)\zeta_{11} =0,\nonumber\\
 \zeta_{12}''   -  
 r_2(x)\zeta_{12} =0,\\
 \zeta_{21}''   -   r_1(x)\zeta_{21} =\tilde{K_2^{(1)}},\nonumber\\
 \zeta_{22}''   -  
 r_2(x)\zeta_{22} =\tilde{ K_2^{(2)}},
\end{eqnarray}
here
\begin{eqnarray*}
  r_1(x) & = &  \frac{4h^6x^6(A+\frac{15}{2})+4A+3}{4x^2(h^6x^6+1)^2},\\
   r_2(x) & = & \frac{9(10h^6x^6+7)}{4x^2(h^6x^6+1)^2}.
\end{eqnarray*}

Now we have

\begin{eqnarray*}
 \tilde{K_2}^{(1)}(\zeta_{11},\zeta_{12})  =   K_2^{(1)}\frac{x^{1/2}}{(h^6x^6+1)^{\frac{1}{4}}} \\
 = \left(4A\frac{\zeta_{11}\zeta_{12}}{x^{3/2}(h^6x^6+1)^{3/4}}+\left(\frac{3D}{2}\right)\frac{(\zeta_{11})^2}{x^{3/2}(h^6x^6+1)^{3/4}}\right)\\
 \tilde{K_2}^{(2)}(\zeta_{11},\zeta_{12})  =   K_2^{(2)}\frac{x^{1/2}}{(h^6x^6+1)^{\frac{1}{4}}}\\
  =  \left(2A\frac{(\zeta_{11})^2}{x^{3/2}(h^6x^6+1)^{3/4}}+30\frac{(\zeta_{12})^2}{x^{3/2}(h^6x^6+1)^{3/4}}\right).
\end{eqnarray*}

Without loss of generality, we can assume that

$\zeta_{11}^{(1)}(\zeta_{11}^{(2)})'-\zeta_{11}^{(2)}(\zeta_{11}^{(1)})'=1$ and $\zeta_{12}^{(1)}(\zeta_{12}^{(2)})'-\zeta_{12}^{(2)}(\zeta_{12}^{(1)})'=1$ . Then the fundamental matrix of (\ref{VE2_zeta}) and its inverse are
\begin{equation}
\label{X}
X (z) =
\begin{pmatrix}
 \zeta_{11} ^{(1)} & \zeta_{11} ^{(2)} & 0 & 0\\
 (\zeta_{11} ^{(1)})' & (\zeta_{11} ^{(2)})'   & 0 & 0 \\
0 & 0 & \zeta_{12} ^{(1)} & \zeta_{12} ^{(2)}\\
0 & 0 & (\zeta_{12} ^{(1)})' & (\zeta_{12} ^{(2)})'
  \end{pmatrix},
  \end{equation}

\begin{equation}
\label{X^1}
X ^{-1}(z) =
\begin{pmatrix}
 (\zeta_{11} ^{(2)})' & -\zeta_{11} ^{(2)} & 0 & 0\\
 -(\zeta_{11} ^{(1)})' & \zeta_{11} ^{(1)}   & 0 & 0 \\
0 & 0 &(\zeta_{12} ^{(2)})' & -\zeta_{12} ^{(2)}\\
0 & 0 & -(\zeta_{12} ^{(1)})' & \zeta_{12} ^{(2)}
  \end{pmatrix}.
  \end{equation}
We show that  a logarithmic term
appears in  local solution of (${{VE}}_2$). For this purpose,
it is sufficient to show that  at least one component of $X^{-1} f_2$
has a nonzero residue at $0$.The appearance of a logarithmic term in a local solution implies that the identity component of the differential Galois group is not Abelian. We calculate
 of $X^{-1} f_2$, which looks like
 $$(-\zeta_{11}^{(2)}\tilde{K_2^{(1)}},\,\zeta_{11}^{(1)}\tilde{K_2^{(1)}},\, -\zeta_{12}^{(2)}\tilde{K_2^{(2)}},\, -\zeta_{12}^{(1)}\tilde{K_2^{(2)}} )^T.$$
Since we have
\begin{equation*}
e^{2\int{\frac{2h^6x^6-1}{x(h^6x^6+1)}}dx}=\frac{(h^6x^6+1)^{1/2}}{x^2},
\end{equation*}
in the above expression, we can replace $\zeta \rightarrow \xi$  and $\tilde{K}_2^{(i)}\rightarrow K_2^{(i)}$, $i=1,\, 2$ and investigate when in its expansion near 0 we have a non-zero term of the second power of $x$. In the expressions for 
$$(-\xi_{11}^{(2)}{K_2^{(1)}},\,\xi_{11}^{(1)}{K_2^{(1)}},\, -\xi_{12}^{(2)}{K_2^{(2)}},\, -\xi_{12}^{(1)}{K_2^{(2)}} )^T,$$
 we look for non-zero coefficients in front of $x^2$.

Let us focus on the first case $\sqrt{A+1}=3k+1$, for $k\in\mathbb{Z}$, then from (\ref{xi_Sol_0})  we get the expansions
\begin{eqnarray}
\label{xi_Case1}
\xi_{11}^{(1)}(x)=x^{-3k}\left( 1+O(x^6)\right),\nonumber\\
\xi_{11}^{(2)}(x)=x^{3k+2}\left( 1+O(x^6)\right),\\
\xi_{12}^{(1)}(x)=x^5\left( 1+O(x^6)\right)\nonumber,\\
\xi_{12}^{(2)}(x)=x^{-3}\left( -20331280+O(x^6)\right).
\end{eqnarray}
We have a non-zero coefficient in front of $x^2$ exactly when $\xi_{11}^{(2)}{K_2^{(1)}}(\xi_{11}^{(1)},\,\xi_{12}^{(1)})$, or  $\xi_{11}^{(2)}{K_2^{(1)}}(\xi_{11}^{(1)},\,\xi_{12}^{(2)})$ and this happens when $k\in 2\mathbb{Z}_{\ge 0}$, and $D\ne 0$.

For the second case $\sqrt{A+1}=-3k+2$, for $k\in\mathbb{Z}$ we have
\begin{eqnarray}
\label{xi_Case2}
\xi_{11}^{(1)}(x)=x^{3k-1}\left( 1+O(x^6)\right),\nonumber\\
\xi_{11}^{(2)}(x)=x^{-3k+3}\left( 1+O(x^6)\right),\\
\xi_{12}^{(1)}(x)=x^5\left( 1+O(x^6)\right)\nonumber,\\
\xi_{12}^{(2)}(x)=x^{-3}\left( -20331280+O(x^6)\right).
\end{eqnarray}
And we have a non-zero coefficient in front of $x^2$ for the term $\xi_{11}^{(1)}{K_2^{(1)}}(\xi_{11}^{(2)},\,\xi_{12}^{(1)})$, whence $k\in 2\mathbb{Z}_{\ge0}+1$ and $D\ne 0$.

For the last case  $\sqrt{A+1}=6k+3$  ($k\in\mathbb{Z}$) the solutions around $0$ are

\begin{eqnarray}
\label{xi_Case3}
\xi_{11}^{(1)}(x)=x^{-6k-2}\left( 1+O(x^6)\right),\nonumber\\
\xi_{11}^{(2)}(x)=x^{6k+4}\left( 1+O(x^6)\right),\\
\xi_{12}^{(1)}(x)=x^5\left( 1+O(x^6)\right)\nonumber,\\
\xi_{12}^{(2)}(x)=x^{-3}\left( -20331280+O(x^6)\right).
\end{eqnarray}
 We obtain a non-zero coefficient in front of $x^2$ for the term $\xi_{11}^{(1)}{K_2^{(1)}}(\xi_{11}^{(1)},\,\xi_{12}^{(1)})$, when $k\in 3\mathbb{Z}_{\ge0}$.

We proved the main

\begin{theorem}
\label{Th_V6}
For $A\ne 3$, the system (\ref{1.2}) does not have an additional meromorphic first integral, if at least one of the following conditions holds:

1. $\sqrt{A+1}\notin\mathbb{Q}$;

 2. $\sqrt{A+1}\in\mathbb{Q}\setminus\{ 3k+1, -3k+2, 6k+3\}$,  for $k\in\mathbb{Z}$;
 
 3. $\sqrt{A+1}\in\{3k+1, -3k+2, 6k+3\}$, for $k\in\mathbb{Z}$;

3.1. In the case, when  $\sqrt{A+1}=3k+1$,  the conditions are:  $k\in 2\mathbb{Z}_{\ge 0}$, and $D\ne 0$;

3.2. In the case, when  $\sqrt{A+1}=-3k+2$, the conditions are:  $k\in 2\mathbb{Z}_{\ge0}+1$, and $D\ne 0$;

3.3. In the case $\sqrt{A+1}=6k+3$,  the conditions are:  $k\in 3\mathbb{Z}_{\ge0}$.
\end{theorem}

\section{Comments}
In fact, the proof is not quite complete and therefore some clarifications are necessary. We need to clarify what happens when the coefficients of (\ref{xi_Sol_0}) vanish or are undefined. This happens exactly when $A=0$, $A=3$ and $A=8$. 

Let us start with the last case $A=8$. In this case, the solutions are expressed as hypergeometric functions with exponents $a=\ds{\frac{2}{3}}$, $b=\ds{\frac{5}{6}}$, and $c=\ds{\frac{1}{2}}$. We apply Kimura's result \cite{Kimura} and obtain non-integrability in this case ($\hat{\nu}=-1$).

Let us focus on case $A=0$. Then, the first variational equation VE1 (only in it there is a difference) is:
\begin{equation}
\label{VE1_inftyA=0}
 \xi_{11}''  + \frac{2h^6x^6-1}{x(h^6x^6+1)}\xi_{11}' =0, 
\end{equation}
with solutions
\begin{eqnarray}
\label{xi_CaseA=0}
\xi_{11}^{(1)}(x)=x^{2}\left( 1+O(x^6)\right),\nonumber\\
\xi_{11}^{(2)}(x)=c\left( 1+O(x^6)\right).
\end{eqnarray}
We notice that the coefficient in front of $x^2$ in the expression $\xi_{11}^{(1)}(\xi_{11}^{(2)})^2$ is different from 0 at $D\ne0$. This condition, as we noted above, is sufficient ground for non-integrability.

Let us consider several well-known integrable and non-integrable cases
and examine how they fit into the results obtained above.
In \cite{Yarmo2} a previously conjectured integrable homogeneous potential of degree~6
was found.
Here, we consider the potential
$$
V_6=r^6+80r^2z^4+24r^4z^2+64z^6.
$$
It is easy to verify that the non-integrability conditions of
Theorem~(\ref{Th_V6}) are not satisfied in this case.
Moreover, the Morales--Ramis necessary conditions for integrability
are also fulfilled.
The Poincar\'e map also exhibits the typical behavior of an integrable system.
Despite these indications of the absence of chaos,
the system turns out to be non-integrable.

Consider the following canonical transformation:
\begin{eqnarray*}
x=r,\, p_x=p_r\\
y=2z,\, p_y=\frac{p_z}{2}.
\end{eqnarray*}
We obtain a Hamiltonian system with potential $V_{6}=x^6+6x^4y^2+5x^2y^4+y^6$, which is found to be non-integrable according to both the results of the present paper and the Morales-Ruiz criterion.  We report a "counterintuitive finding": although all criteria suggested the existence of an additional first integral for the Hamiltonian system, the system was ultimately shown to be "non-integrable". However, this conclusion does not hold for the following reason: throughout our analysis, we have considered only Hamiltonian systems of the form 
$$H=\frac12\left(p_r^2+p_z^2\right)+V(r,z).$$
Following the canonical transformation described above, the system is transformed into the form $$H=\frac{1}{2}(p_x^{2}+2p_y^{2})
+x^6+6x^4y^2+5x^2y^4+y^6 . $$
Since this system does not belong to the class under consideration, the results established in this paper, as well as the Morales–Ramis method, cannot be applied to it.

Moreover, the transformed system admits an additional independent first integral given by 
$$I=\frac{y}{4}\,p_x^{2}
-x\,p_xp_y
-\frac{3}{2}x^6y
-2x^4y^3
-\frac{1}{2}x^2y^5 . $$
It is straightforward to verify $\{H,I\}=0$.

\begin{figure}[ht]
    \centering
    \begin{subfigure}[b]{0.45\textwidth}
        \centering
        \includegraphics[width=\textwidth]{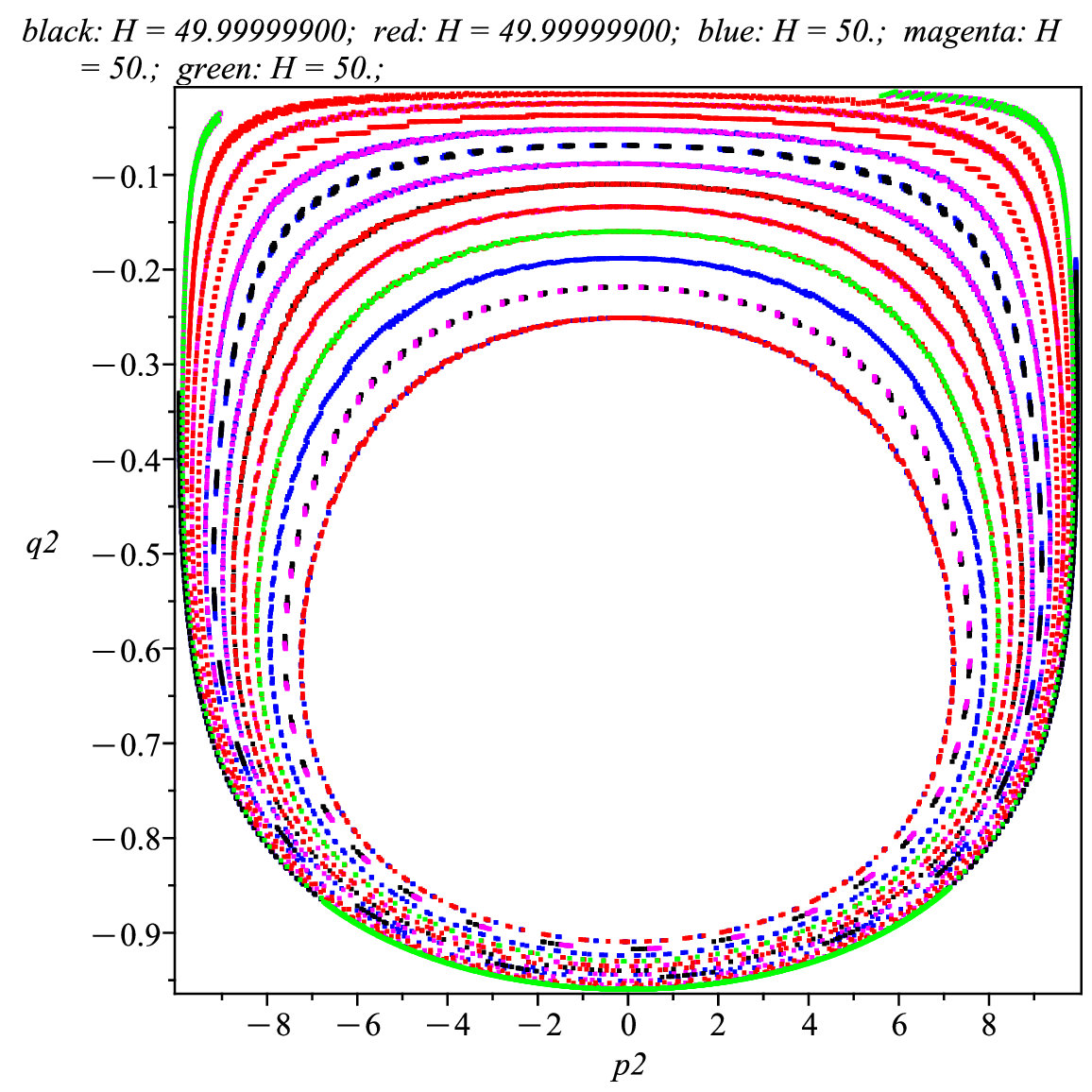}
        \caption{$(p_r,r)$ - projection;}
    \end{subfigure}
    \hfill 
    \begin{subfigure}[b]{0.45\textwidth}
        \centering
        \includegraphics[width=\textwidth]{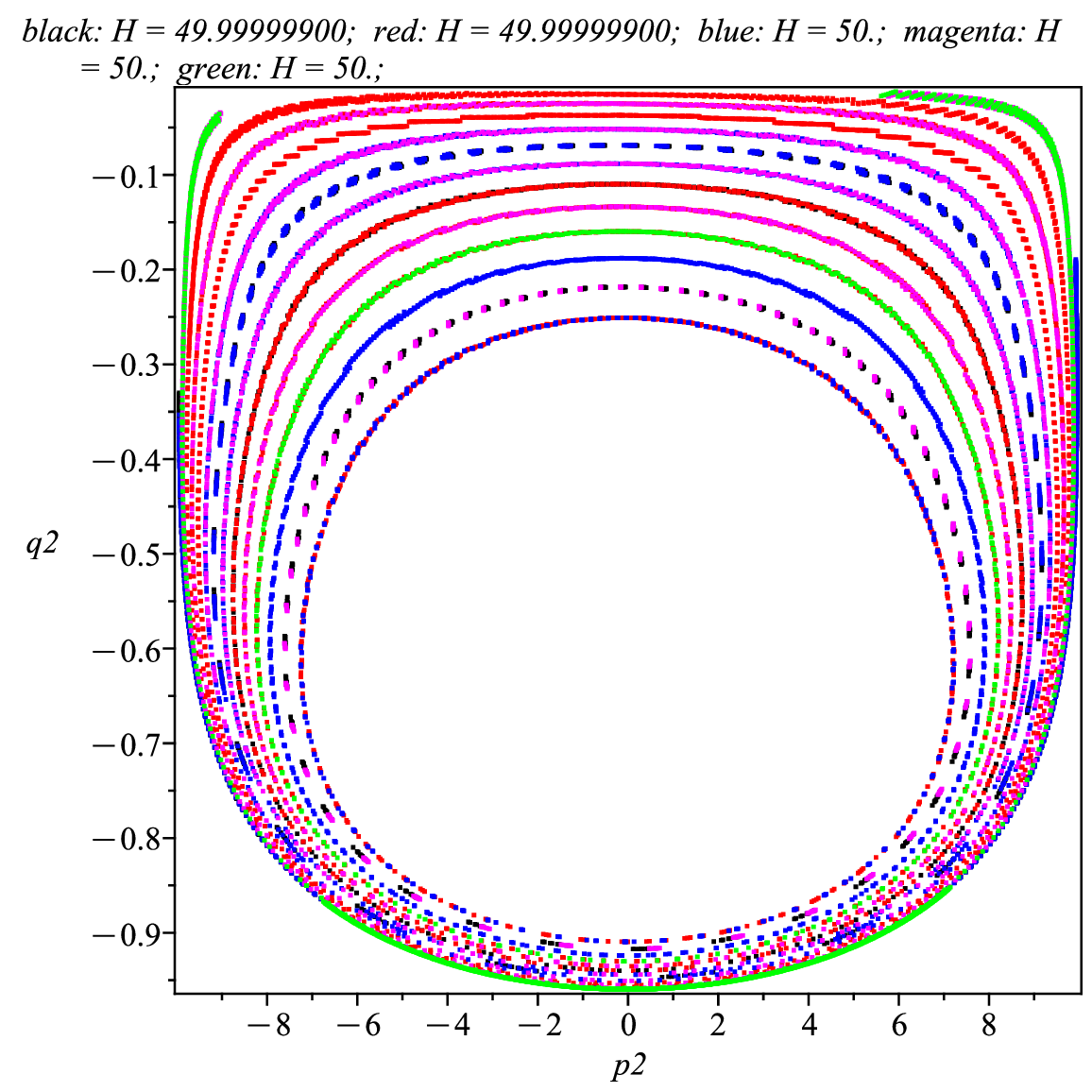}
        \caption{$(p_z,z)$ - projection.}
    \end{subfigure}
\caption{Poincare--section in the case $V_6=r^6+80r^2 z^4+24r^4 z^2+64z^6$ .}
    \label{fig:my_plots6}
\end{figure}

We will now examine several cases that are resolved by finding the eigenvalues of the Hessian, and compare them with the results obtained from the Theorem(\ref{Th_V6}) we have already established.

Let us briefly describe the method for investigating the integrability of homogeneous potentials using Darboux points (Morales-Ruiz-method). (See \cite{Yoshida0}, \cite{Yoshida1}, \cite{MR1.5} and \cite{MPr1} for details.) We restrict our considerations to two degrees of freedom and homogeneous potential of degree 6. Point $c\in \mathbb{C}^2\setminus\{0\}$ is called the Darboux point for the homogeneous potential $V(q)$, $q\in \mathbb{C}^2$ if $V'(c)=c$ is satisfied. The Hamiltonian system with such a potential is integrable with meromorphic first integrals if the second eigenvalue  of the Hessian $V''(c)$, ($\lambda_1=5$), belongs to one of the families $\{p+3p(p-1)$ or  $\frac{1}{2}(\frac{5}{6}+6p(p+1))\}$,   for $p\in\mathbb{Z}$.

Let us look at a potential that is not very difficult to study $V=r^3z^3$. For it we have $A=0$ and by Theorem(\ref{Th_V6}), it is non-integrable. Let us find the Darboux points in this case. We have 
the points $(r,\, z)=(3^{-1/4},\,3^{-1/4}) $, $(r,\, z)=(-(-3)^{-1/4},\,(-3)^{-1/4}) $ and  $(r,\, z)=((-3)^{-1/4},\,-(-3)^{-1/4})$.
The eigenvalues of the Hessian at the Darboux points are $\lambda_1=5$ and $\lambda_2= -1$. We verify immediately that
$$-1\notin \{p+3p(p-1),\, \frac{1}{2}(\frac{5}{6}+6p(p+1))\},\, p\in\mathbb{Z}.$$
Hence, by the Morales–Ramis criterion, the potential is non-integrable.

We now proceed to discuss the case $A=3$. We start with an example:\\
$V=(r^2+z^2)^3$. This potential is excluded from the results of Theorem(\ref{Th_V6}), and it also has no Darboux points. This represents a typical borderline case. The function $L=zp_r-rp_z$ is the first integral for the system with this potential. This example demonstrates, as noted in \cite{Yarmo2} and \cite{Valls}, that any potential of the form $V=f(r^2+z^2)$ is integrable. It is easy to check that the function $L=zp_r-rp_z$ is the first integral for the system with a potential $V=f(r^2+z^2)$, if the function $f$ has a derivative, of course.  

Let us show another example of $A=3$.

 Let $V=(r^2+z^2)^3+Dr^3z^3+Erz^5$, here we assume that $(D,\, E)\ne (0,\, 0)$. We will show that a large part of the parameters $\alpha=\frac{D}{E}$, the potential is non-integrable. The first eigenvalue of the Hessian $V''(c)$ is the standard 5, and for the second we have two possibilities $\lambda_{\pm}(\alpha )=\frac{(21\alpha+25\mp 5\sqrt{9\alpha^2-18\alpha+25})}{6\alpha}$. If we assume $\sqrt{9\alpha^2-18\alpha+25}\notin \mathbb{Q}$, then we clearly have non-integrability. When $\sqrt{9\alpha^2-18\alpha+25}\in \mathbb{Q}$, we have $\alpha=\frac{16-k^2}{6k}+1$, for $k\in \mathbb{Q}$, we get for
$$\lambda_{+}(\alpha)=\lambda_{+}(k)=6-\frac{10}{k+2}$$
 and 
$$\lambda_{-}(\alpha)=\lambda_{-}(k)=1-\frac{40}{k-8}.$$
We apply the Morales--Ruiz--Ramis table and obtain the conditions for non-integrability in this case:
\begin{theorem}
\label{Th_alphaDE}

For $\alpha =\frac{D}{E}=\frac{16-k^2}{6k}+1$,  and $p_i\in \mathbb{Z}$,  $i=1,\, 2, \,3,\, 4$ a necessary condition for the non-integrability of potential $V=(r^2+z^2)^3+Dr^3z^3+Erz^5$ is that at least one of the conditions is satisfied;

1. $k\ne \frac{10}{6-4p_1-3p_1^2}-2$;

2. $k\ne \frac{40}{1-4p_2-3p_2^2}+8$;

3. $k\ne \frac{120}{67-36p_3-36p_3^2}-2$;

4. $k\ne \frac{480}{7-36p_4-36p_4^2}+8$.

\end{theorem}

The necessary condition for integrability is that both $\lambda_{+}(k)$ and $\lambda_{-}(k)$ ($k\in\mathbb{Q}$) are in one of the families $M_p:=\{p+3p(p-1)\}$ or $K_p:=\{\frac{1}{2}(\frac{5}{6}+6p(p+1))\}$, for   $p\in \mathbb{Z}$. (Here we assume a fixed $k\in\mathbb{Q}$ and a variable $p\in \mathbb{Z}$.) However, this is impossible for $k\in\mathbb{Q}$ and $p\in \mathbb{Z}$. 
It turns out that since the negations of these conditions (1., 2., 3. and 4.) are not compatible, we can conclude that the potential $V=(r^2+z^2)^3+Dr^3z^3+Erz^5$  is non-integrable in the case under consideration.
The same is shown by the Poincare section for  $D=1$, $E=1$, and $D=2$, $E=1$. 
\begin{figure}[ht]
    \centering
    \begin{subfigure}[b]{0.45\textwidth}
        \centering
        \includegraphics[width=\textwidth]{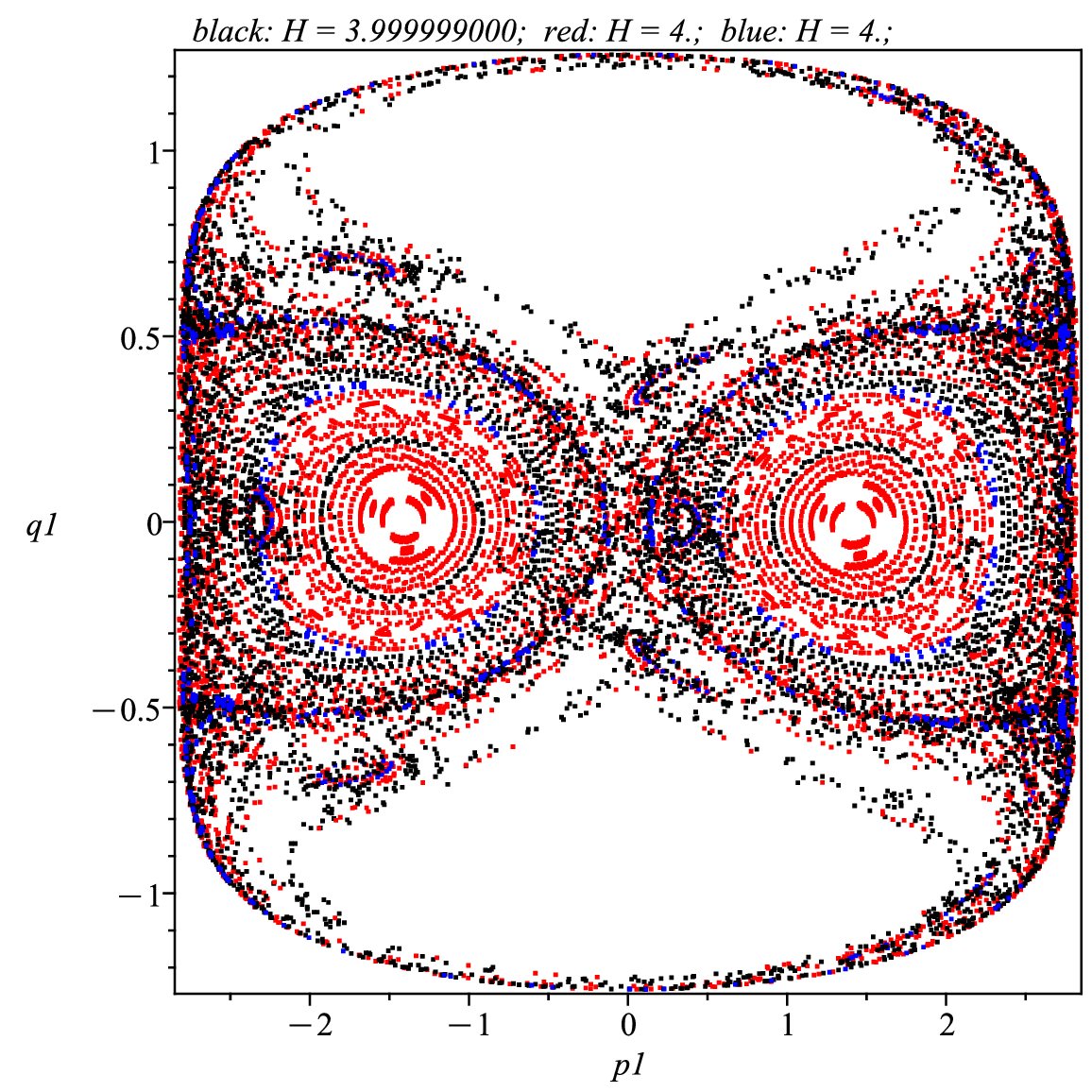}
        \caption{$(p_r,r)$ - projection;}
    \end{subfigure}
    \hfill 
    \begin{subfigure}[b]{0.45\textwidth}
        \centering
        \includegraphics[width=\textwidth]{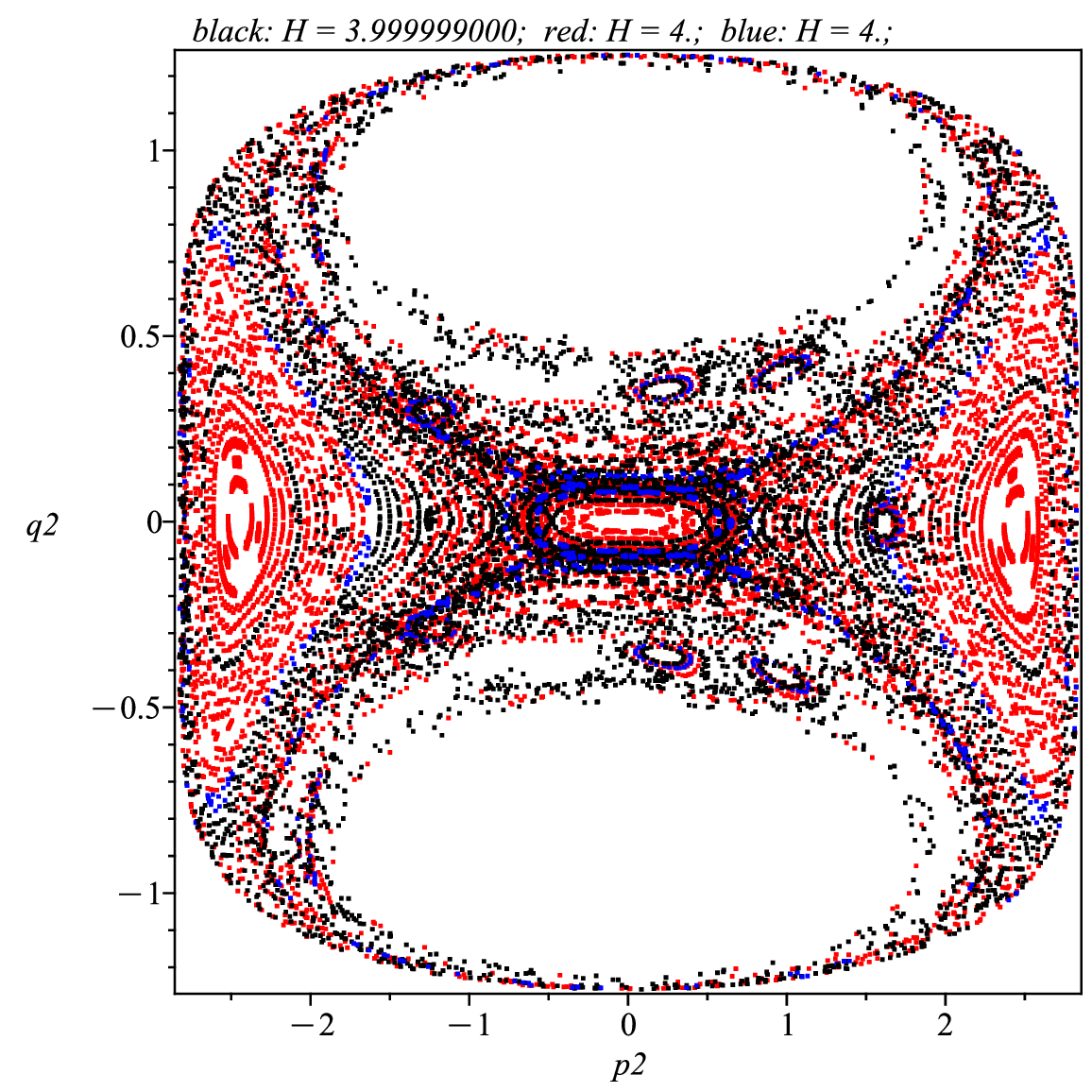}
        \caption{$(p_z,z)$ - projection.}
    \end{subfigure}
    
    \caption{Poincare--section in the case $D=1$ and $E=1$.}
    \label{fig:my_plots1}
\end{figure}
\begin{figure}[ht]
    \centering
    \begin{subfigure}[b]{0.45\textwidth}
        \centering
        \includegraphics[width=\textwidth]{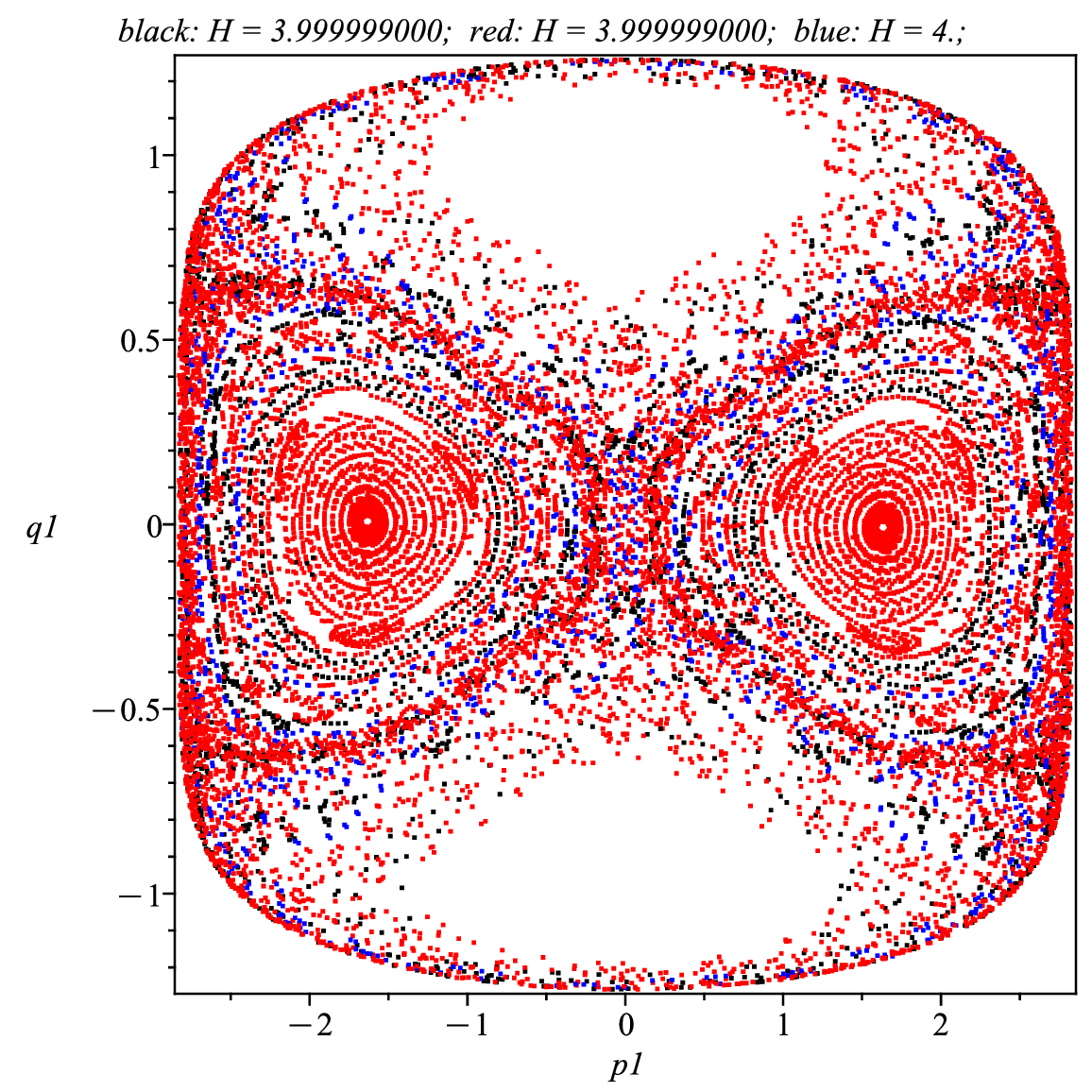}
        \caption{$(p_r,r)$ - projection;}
    \end{subfigure}
    \hfill 
    \begin{subfigure}[b]{0.45\textwidth}
        \centering
        \includegraphics[width=\textwidth]{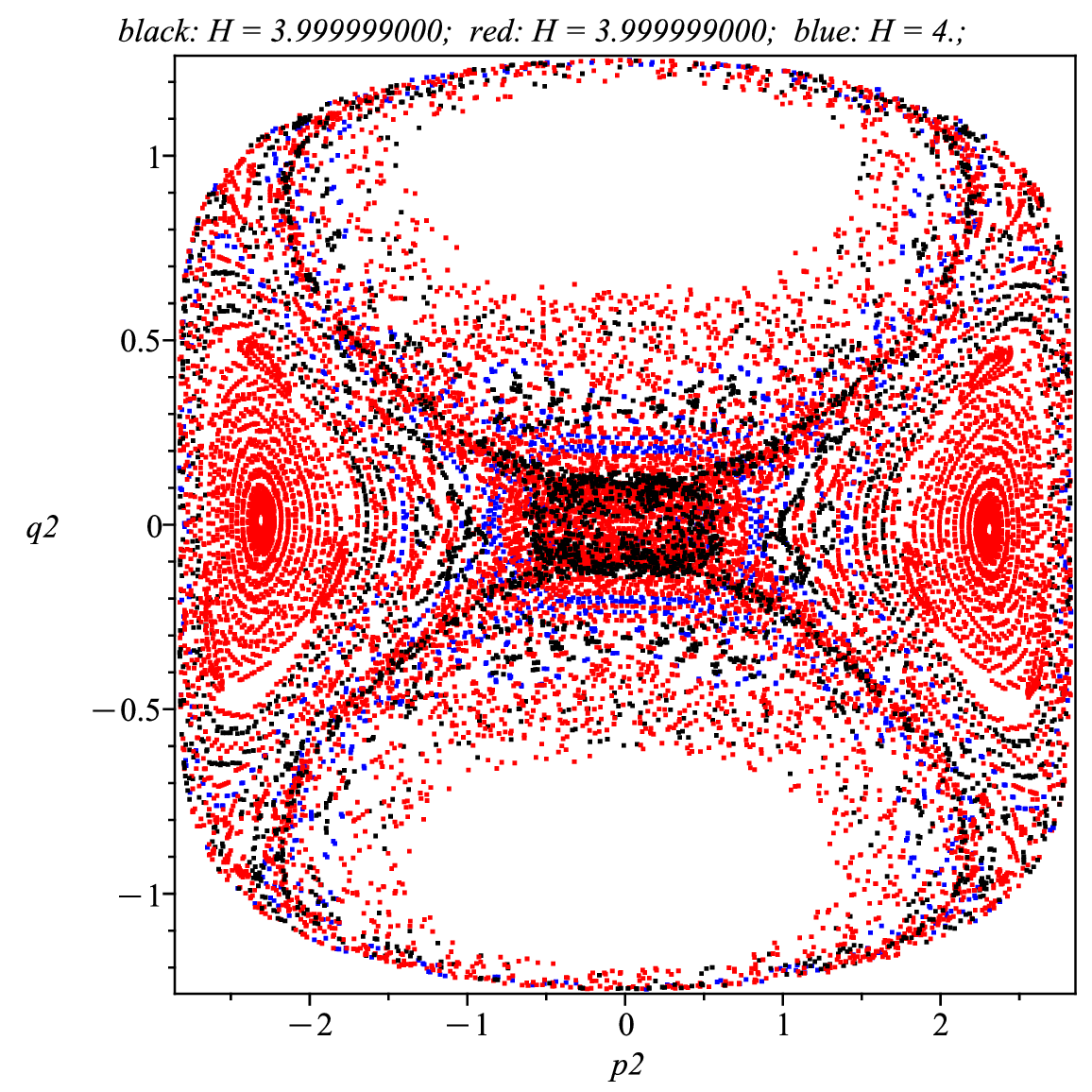}
        \caption{$(p_z,z)$ - projection.}
    \end{subfigure}
    
    \caption{Poincare--section in the case $D=2$ and $E=1$.}
    \label{fig:my_plots2}
\end{figure}

\begin{figure}[ht]
    \centering
    \begin{subfigure}[b]{0.45\textwidth}
        \centering
        \includegraphics[width=\textwidth]{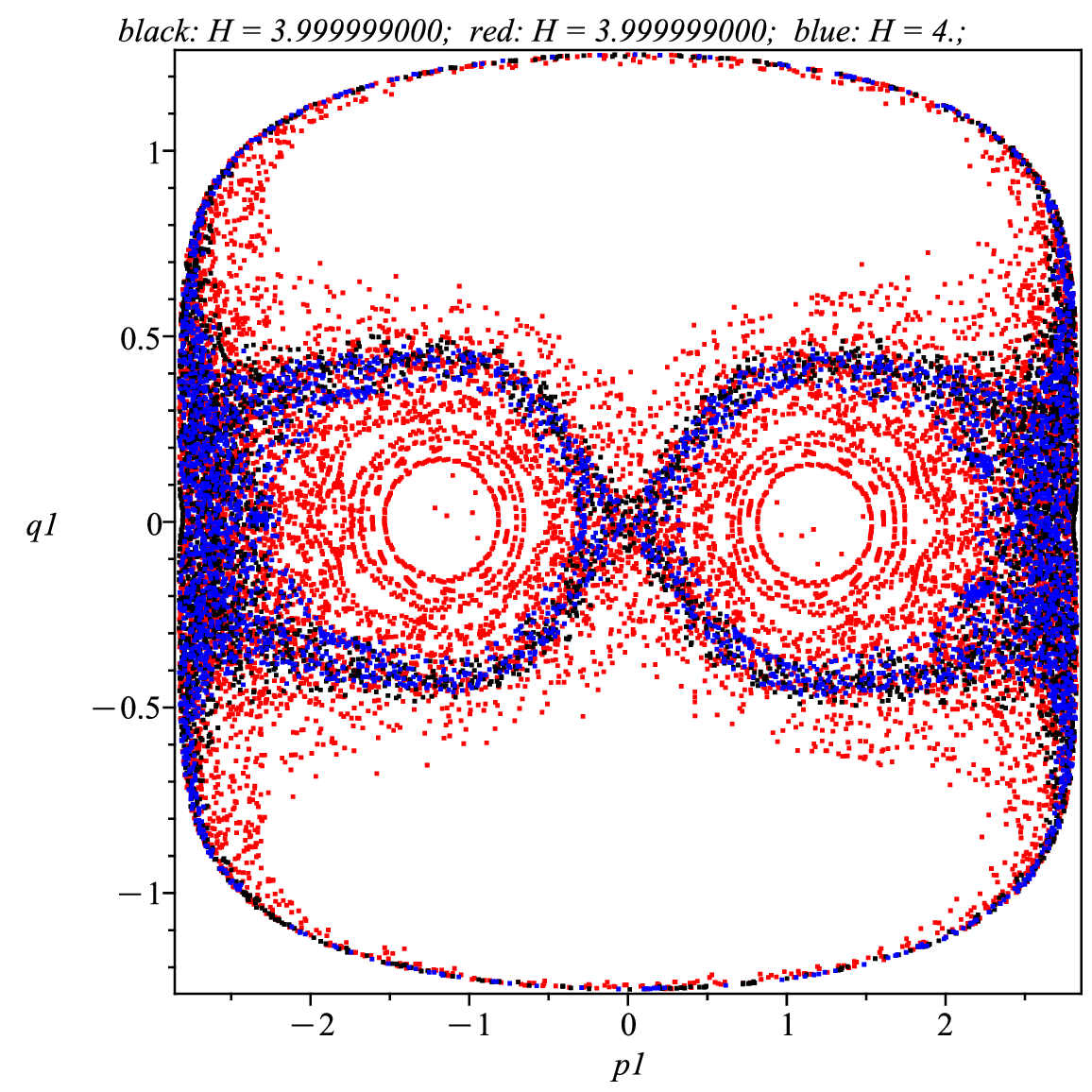}
        \caption{$(p_r,r)$ - projection;}
    \end{subfigure}
    \hfill 
    \begin{subfigure}[b]{0.45\textwidth}
        \centering
        \includegraphics[width=\textwidth]{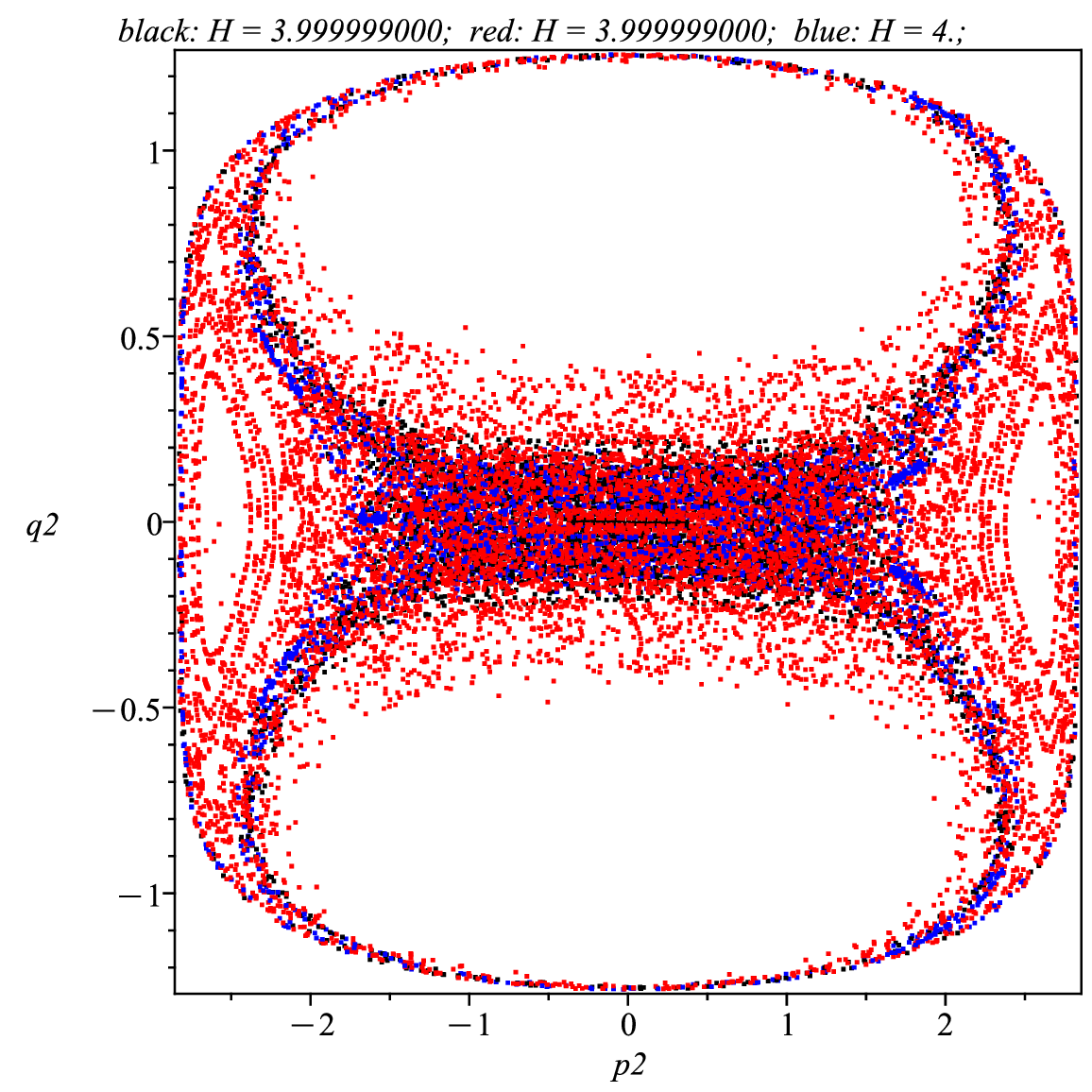}
        \caption{$(p_z,z)$ - projection.}
    \end{subfigure}
    
    \caption{Poincare--section in the
case $D=0$ and $E=1$.}
    \label{fig:my_plots3}
\end{figure}

At $D=0$ the Poincare section also exhibits chaotic dynamics.

Let us consider another example to apply the established results regarding the integrability of homogeneous potentials of degree 6. This example serves as the fundamental motivation for our investigation into this task.  We examine homogeneous potential $V_6 =r^6+\mu r^4z^2+\delta r^2z^4+z^6$ with a particular focus on the necessary conditions for its integrability. These potentials arise in the study of double nonlinear oscillators. According to Theorem (\ref{Th_V6}), conditions $\mu=k^2-1$ and $\delta=m^2-1$ must be satisfied for $k,\, m \in \mathbb{Z}$ to hold. Applying the Morales-Ruiz criterion leads us to $|k|=|m|\in\{1, \, 2,\, 4\}$. Hence, we obtain that the only possible integrable potentials of the considered form are: 

1.  $V_{6_1}=r^6+z^6$ - here the variables separate, i.e. the system is integrable;

2. $V_{6_2}=r^6+3r^4z^2+3r^2z^4+z^6=(r^2+z^2)^3$ - this case was considered above and also turns out to be integrable;

3. $V_{6_3 }=r^6+15r^4z^2+15r^2z^4+z^6$ - here the situation is somewhat more complicated, but the canonical transformation 
\begin{eqnarray*}
q_1=\frac{r+z}{\sqrt{2}},\, p_1=\frac{p_r+p_z}{\sqrt{2}},\\
q_2=\frac{r-z}{\sqrt{2}},\, p_2=\frac{p_r-p_z}{\sqrt{2}},
\end{eqnarray*}
 reduces the problem to a  system with Hamiltonian
$$H=\frac{1}{2}(p_1^2+p_2^2)+4q_1^6+4q_2^6, $$
the system is integrable.

We have thus obtained that these are the only integrable potentials of the form  $V_6 =r^6+\mu r^4z^2+\delta r^2z^4+z^6$.

The use of both approaches provides a rapid answer to the question of the necessary and sufficient conditions for the meromorphic integrability of the sixth-degree Henon–Heiles model $V_6 =r^6+\kappa r^3z^3+z^6$. This system is integrable if and only if $\kappa=0$.

 The case $A=3$, unfortunately, remained not fully clarified, but our attempts to solve  it with variations higher than the third did not yield conclusive results.

If we use the homogeneity of the potential $V_6$, and we perform the transformation $z\leftrightarrow r$, and shift $A\leftrightarrow \frac{B}{C}$, $C\ne 0$ we may get a similar, a little more general result. 
Since no essential part of the argument relies on $A$, $B$, $C$, $D$ and $E$ being real constants, they can also be assumed to be complex numbers.

\section{Acknowledgements}
 This work was partially supported by Sofia University under Contract No. 80-10-32/01.04.2026.
\vspace{5ex}
\FloatBarrier
%\clearpage
\appendix

 \section{Differential Galois Theory}

In the following, we review several classical and advanced results that will be used in our subsequent analysis.

We say that a system with $n$ degrees of freedom is integrable in the sense of Liouville
if it possesses a complete set of $n$ independent first integrals
$f_1=H,\, f_2,\, \dots,\, f_n$ in involution, i.e.,
the Poisson brackets $\{f_i,\,f_j\}=0$, for $i\ne j$

 We briefly recall several facts concerning the integrability of Hamiltonian systems
in the complex domain, the Ziglin--Morales--Ramis theory,
and its relation to the differential Galois groups of linear differential equations.
Our presentation follows \cite{MR1,ChG,MR2,MRS1}..

We assume that a Hamiltonian system
\begin{equation}
\label{2.1}
\dot{x} = X_{H} (x), \quad t \in \mathbb{C}, \quad x \in M
\end{equation}
corresponds to an analytic Hamiltonian $H$, defined on the complex
$2 n$-dimensional manifold $M$. Suppose that system~(\ref{2.1}) admits a non-equilibrium solution
$\Psi(t)$, and denote by $\Gamma$ its phase curve.
The variational equation (VE) along $\Psi(t)$ takes the form
\begin{equation}
\label{2.2}
\dot{\mathbf{\xi}} = D X_{H} ( \Psi (t)) \mathbf{\xi},
\quad \mathbf{\xi} \in T_{\Gamma} M.
\end{equation}

Next, we consider the normal bundle
$F:=T_{\Gamma}M/T\Gamma$ of $\Gamma$
and let $\pi:T_{\Gamma}M\to F$ denote the natural projection.
Equation~(\ref{2.2}) induces a differential equation on $F$
\begin{equation}
\label{2.3}
\dot{\eta} = \pi_{*} (D X_{H} ( \Psi (t))(\pi^{-1} \eta) , \quad \eta \in F.
\end{equation}
which is  called a normal variational equation (NVE) around $\Gamma$.

The normal variational equation (NVE)~(\ref{2.3})
possesses a first integral $dH$ linear on the fibres of $F$.
For each $r\in\mathbb{C}$, the level set
$
F_r:=\{\eta\in F \mid dH(\eta)=r\}
$
forms a $(2n-2)$-dimensional affine bundle over $\Gamma$.

%%%%%%
We will call $F_r$ the reduced phase space of (\ref{2.3}) and the restriction of
the (NVE) on $F_r$ is called the reduced normal variational equation.

The main result of the Morales--Ramis theory \cite{MR1} is the following:
\begin{theorem}
\label{th2}
Let us assume that the Hamiltonian system (\ref{2.1}) has $n$ meromorphic first integrals in involution,
then the identity component $G^0$ of the Galois group of the variational equation is abelian.
\end{theorem}

Next we consider a linear  system
\begin{equation}
\label{2.10}
y' = A (x) y, \quad y \in \mathbb{C}^n ,
\end{equation}
or  linear homogeneous differential equation, which is essentially the same
\begin{equation}
\label{2.11}
y^{(n)} + a_1 (x) y^{(n-1)} + \ldots + a_n (x) y = 0,
\end{equation}
with $x \in \mathbb{CP}^1$  and $A \in \mathrm{gl} (n, \mathbb{C} (x))$,
($a_j (x) \in \mathbb{C} (x))$.
Let $S:=\{x_1, \ldots, x_s\}$ be the set of singular points of (\ref{2.10}) (or (\ref{2.11})) and let $Y (x)$ be a fundamental solution
of (\ref{2.10})  (or (\ref{2.11}))  at $x_0 \in \mathbb{C} \setminus S$. By the existence theorem, this solution is analytic
in a neighbourhood of $x_0$.
The analytic continuation of $Y(x)$
along a nontrivial loop in $\mathbb{CP}^1$
defines a linear automorphism of the space of solutions,
called the monodromy operator. Analytically, this transformation is described as follows.
To each loop
$\gamma \in \pi_1(\mathbb{CP}^1 \setminus S, x_0)$
there corresponds a linear automorphism $\Delta_{\gamma}$,
which acts on $Y(x)$ by right multiplication with a constant matrix
$M_{\gamma}$, called the monodromy matrix
$$
\Delta_{\gamma} Y (x) = Y (x) M_{\gamma}.
$$
The collection of all such matrices forms the monodromy group.

We add another object to the  (\ref{2.10}) (or (\ref{2.11})) - a differential
Galois group.
We have a differential field $K$, that is a field with a derivation $\partial = '$, i.e.
an additive mapping complying with the derivation rule. Differential automorphism
of $K$ is an auto\-mor\-phism commuting with the derivation.

The coefficient field in (\ref{2.10}) (and  (\ref{2.11})) is $K = \mathbb{C} (x)$. Let $y_{i j}$ be
elements of the fundamental matrix $Y (x)$. Let $L (y_{i j})$ be the
extension of $K$ generated by $K$ and $y_{i j}$ -- a differential
field. This extension is called a Picard--Vessiot's extension.
 Similarly to classical Galois Theory we define the Galois group
$G := Gal_{K} (L) = Gal (L/K)$ to be the group of all differential
automorphisms of $L$ leaving the elements of $K$ fixed.
 Galois group is an algebraic group. It has an unique
connected component $G^0$ which contains the identity and  is a normal
subgroup of finite index.  Galois group $G$ can be represented
as an algebraic linear subgroup of $\mathrm{GL} (n, \mathbb{C})$ by
$$
\sigma (Y (x)) = Y (x) R_{\sigma},
$$
where $\sigma \in G$ and $R_{\sigma} \in \mathrm{GL} (n, \mathbb{C})$.

We can do the same locally at $a \in \mathbb{CP}^1$, replacing $\mathbb{C} (x)$ by the field of germs of meromorphic
functions at $a$. Hence we can speak of a local differential Galois group $G_a$ of (\ref{2.10}) at $a \in \mathbb{CP}^1$,
defined in the same way for Picard-Vessiot extensions of the field $\mathbb{C} \{x-a\}[(x-a)^{-1}]$.

It is worth here that by its definition the monodromy group is contained
in the differential Galois group of the corresponding system.

Next, we are presenting some facts from the theory of linear systems with singularities.
We call a singular point $x_i$  regular if any of the solutions of (\ref{2.10})
(or of (\ref{2.11})) has at most polynomial growth in arbitrary sector with a vertex at
$x_i$. Otherwise the singular point is called  irregular.

We say that the system (\ref{2.10}) has a singularity of the Fuchs type at $x_i$ if $A (x)$ has a
simple pole at $x = x_i$. For the equation (\ref{2.11}) the Fuchs type singularity at  $x_i$
means that the functions $(x - x_i)^j a_j (x)$ are holomorphic in a neighborhood of $x_i$.

If the system (\ref{2.10}) has a singularity of the Fuchs type, then  this singularity is regular.
The opposite is not true. However, for the equation (\ref{2.11}) the regular singularities coincide with
the singularities of Fuchs type.

A system with only regular singularities is called Fuchsian system.
For such systems we have :
\begin{theorem}
\label{th3}
({\bf Schlesinger} )
The differential Galois group coincides with the  Zariski closure in $\mathrm{GL}(n,\mathbb{C})$ of the mo\-no\-dro\-my
group.
\end{theorem}

The fact that $G^0$ is abelian  doesn't imply  necessarily
integrability of the Hamiltonian system.
There is a method which, in the case of abelian Galois group, can draw conclusion when the system (\ref{2.1}) is non-integrable. This method based on
the higher variational equations has been introduced in \cite{MR1}
and  the Theorem \ref{th2} has been extended in \cite{MRS1}. What is the idea of higher variational
equations? For the system (\ref{2.2}) with a particular solution
$\Psi (t)$ we put
\begin{equation}
\label{2.5}
x = \Psi (t) + \varepsilon \xi^{(1)} + \varepsilon^2
\xi^{(2)} + \ldots + \varepsilon^k \xi^{(k)} + \ldots,
\end{equation}
where $\varepsilon$ is a small   parameter. When substituting the
above expression into eq. (\ref{2.2}) and comparing terms with the
same order in $\varepsilon$ we obtain the following chain of
linear non-homogeneous equations
\begin{equation}
\label{2.6}
\dot{\xi}^{(k)} = A (t) \xi^{(k)} + f_k (\xi^{(1)},
\ldots, \xi^{(k-1)}), \quad k = 1, 2, \ldots ,
\end{equation}
where $A (t) = D X_{H} ( \Psi (t))$ and $f_1 \equiv 0$. The
equation (\ref{2.6}) is called k-th variational equation
(${\rm{VE}}_k$). Let $X (t)$ be the fundamental matrix of
(${\rm{VE}}_1$)
$$
\dot{X} = A (t) X .
$$
Then the solutions of $({\rm{VE}}_k), k > 1$ can be expressed with
\begin{equation}
\label{2.7}
\xi^{(k)} = X (t) c (t),
\end{equation}
where $c (t)$ is a solution of
\begin{equation}
\label{2.8}
\dot{c} = X^{-1} (t) f_k .
\end{equation}
Although (${\rm{VE}}_k$) are not actually homogeneous equations,
they can be placed in this framework, and therefore,  successive
extensions $ K \subset L_1 \subset L_2 \subset \ldots \subset
L_k$ can be outlined, where $L_k$ is the extension obtained by adjoining the
solutions of (${\rm{VE}}_k$). The
differential Galois groups $Gal (L_1 /K), \ldots, Gal (L_k /K)$ can be defined accordingly. The following
result is proven in \cite{MRS1}.
\begin{theorem}
({\bf  Morales-Ruiz, Ramis,  Sim\'o})
\label{th3b}
If the Hamiltonian system (\ref{2.2}) is integrable
in Liouville sense then the identity component of each Galois
group $ Gal (L_k /K)$ is abelian.
\end{theorem}

Pay attention that we apply Theorem \ref{th3} to the situation where the
identity component of the Galois group $Gal(L_1/K)$ is abelian.
This means that the first variational equation is solvable. Once
we have the solution of $({\rm{VE}}_1)$, then the solutions of
$({\rm{VE}}_k)$ can be found by the method of 
constant variations as explained above. Hence, the differential Galois groups $Gal(L_k /K)$ are solvable.
One possible way to show that some of them are not commutative is
to find a logarithmic term in the corresponding solution.
We need to explain why the existence of a non-zero logarithmic term in $VE_k$ around some singular point guarantees us non-integrability. The Galois group $Gal(L_k /K)$ is abelian, if and only if,  the local monodromy of the $({\rm{VE}}_k)$ around the singular point of the coefficients is identity. If for some $k$, we obtain non-zero residue in the Laurent expansions of the expressions of $ X^{-1} (t) f_k$, near singularity point, then the local monodromy will be represented by  a lower (or upper) triangular matrix which is not an identity, i. e. the Galois group $Gal(L_k /K)$ is not abelian (see
detailed descriptions and explanations in \cite{MR1,MR2,MRS1}).

\section*{AUTHOR DECLARATIONS}
\subsection*{Conflict of Interest}
The authors have no conflicts to disclose.
\subsection*{Author Contributions}

\textbf{D. Neykova:}
 Formal analysis (supporting);
Investigation (equal); Methodology (supporting);
Validation (lead); Visualization (equal);
Writing -- review \& editing (equal).

\textbf{G. Georgiev:}
Conceptualization (lead); Formal analysis (lead);
Investigation (equal); Methodology (lead);
Project administration (lead); Supervision (lead);
Writing -- original draft (lead);
Writing -- review \& editing (equal)
\nocite{*}
\section*{References}
\bibliographystyle{aipnum4-1} 
\bibliography{SixPot}% Produces the bibliography via BibTeX.

\end{document}